\documentclass[11pt]{article}
\usepackage{graphicx}
\usepackage{amsfonts}
\usepackage{amsmath, amscd, amsfonts, amssymb, graphicx, color, enumerate}
\usepackage[bookmarksnumbered, colorlinks, plainpages]{hyperref}
\hypersetup{colorlinks=true,linkcolor=red, anchorcolor=green, citecolor=cyan, urlcolor=red, filecolor=magenta, pdftoolbar=true}
\usepackage{theorem}
\newtheorem{Lem}{Lemma}[section]

\newtheorem{theorem}[Lem]{Theorem}
\newtheorem{proposition}[Lem]{Proposition}
\newtheorem{corollary}[Lem]{Corollary}

\newtheorem{remark}[Lem]{Remark}

\newcommand{\qed}{\hbox{\rule{6pt}{6pt}}}

\begin{document}
\title{Further improvements of Young inequality}
\author{Shigeru Furuichi$^1$\footnote{E-mail:furuichi@chs.nihon-u.ac.jp}\\
$^1${\small  Department of Information Science, College of Humanities and Sciences, Nihon University,}\\
{\small 3-25-40, Sakurajyousui, Setagaya-ku, Tokyo, 156-8550, Japan}}
\date{}
\maketitle
{\bf Abstract.}
We focus on the improvements for Young inequality. We give elementary proof for known results by Dragomir, and we give remarkable notes and some comparisons. Finally, we give new inequalities which are extensions and improvements for the inequalities shown by Dragomir.
\vspace{3mm}


\vspace{3mm}
\vspace{3mm}

%%%%%%%%%%%%%%%%%%%%%%%%%%%%%%%%%%%%%%%%%%%%%%%%%%%%%%%%%%%%%%%%%%%%%%%%%%%%%%%%%%%%%%%%%%%%%%%%%%%%%%%%%%%%%%%%%%%%%%%%%%%%%%%%%%%%%%%%%%%%%%%%%%%%%%%%%%%%%%%
%%%%%%%%%%%%%%%%%%%%%%%%%%%%%%%%%%%%%%%%%%%%%%%%%%%%%%%%%%%%%%%%%%%%%%%%%%%%%%%%%%%%%%%%%%%%%%%%%%%%%%%%%%%%%%%%%%%%%%%%%%%%%%%%%%%%%%%%%%%%%%%%%%%%%%%%%%%%%%%
%%%%%%%%%%%%%%%%%%%%%%%%%%%%%%%%%%%%%%%%%%%%%%%%%%%%%%%%%%%%%%%%%%%%%%%%%%%%%%%%%%%%%%%%%%%%%%%%%%%%%%%%%%%%%%%%%%%%%%%%%%%%%%%%%%%%%%%%%%%%%%%%%%%%%%%%%%%%%%%
%%%%%%%%%%%%%%%%%%%%%%%%%%%%%%%%%%%%%%%%%%%  Section1  %%%%%%%%%%%%%%%%%%%%%%%%%%%%%%%%%%%%%%%%%%%%%%%%%%%%%%%%%%%%%%%%%%%%%%%%%%%%%%%%%%%%%%%%%%%%%%%%%%%%%%%%
%%%%%%%%%%%%%%%%%%%%%%%%%%%%%%%%%%%%%%%%%%%%%%%%%%%%%%%%%%%%%%%%%%%%%%%%%%%%%%%%%%%%%%%%%%%%%%%%%%%%%%%%%%%%%%%%%%%%%%%%%%%%%%%%%%%%%%%%%%%%%%%%%%%%%%%%%%%%%%%
%%%%%%%%%%%%%%%%%%%%%%%%%%%%%%%%%%%%%%%%%%%%%%%%%%%%%%%%%%%%%%%%%%%%%%%%%%%%%%%%%%%%%%%%%%%%%%%%%%%%%%%%%%%%%%%%%%%%%%%%%%%%%%%%%%%%%%%%%%%%%%%%%%%%%%%%%%%%%%%
%%%%%%%%%%%%%%%%%%%%%%%%%%%%%%%%%%%%%%%%%%%%%%%%%%%%%%%%%%%%%%%%%%%%%%%%%%%%%%%%%%%%%%%%%%%%%%%%%%%%%%%%%%%%%%%%%%%%%%%%%%%%%%%%%%%%%%%%%%%%%%%%%%%%%%%%%%%%%%%
\section{Introduction}
For $a ,b >0$ and $0 \leq v \leq 1$, the inequality 
$$a^{1-v}b^v \leq (1-v) a +vb$$
 holds and it is called Young inequality. 
This inequality is simplified as
$$t^v \leq (1-v) +vt$$ 
for $t >0$ and  $0 \leq v \leq 1$.  We use this simplified notation for some inequalities throughout this paper. Recently, a number of refinements for Young inequality are studied. In this paper, we focus on the refinements for Young inequality by Dragomir. We give alternative proofs of refined Young inequalities given in \cite{Dragomir01,Dragomir02}, with elementary calculations.
We also show the inequalities we proved in the previous paper \cite{FM}, 
give better estimates than ones proved in \cite{Dragomir02}.
Finally we extend and improve the inequalities given by Dragomir in \cite{Dragomir01,Dragomir02}. We also give the inequalities for the operator version as corollaries.

\section{Some remarks for recent results}
Recently, Dragomir established the following refinement of Young inequality in \cite{Dragomir01}.
\begin{theorem} {\bf (\cite{Dragomir01})}    \label{theorem_2_1}
For $t >0$ and $0 \leq v \leq 1$, 
\begin{equation} \label{ineq01_theorem01}
\frac{(1-v) + vt}{t^v} \leq \exp\left( v(1-v) \frac{(t-1)^2}{t} \right)
\end{equation}
\end{theorem}

{\it Proof}:
To prove the inequality (\ref{ineq01_theorem01}), we put
$$
f_v(t) \equiv (1-v) t^{-v}+v t^{1-v} - \exp\left( v(1-v) \frac{(t-1)^2}{t} \right).
$$
Then we calculate
$$
\frac{df_v(t)}{dt} = v(1-v)(1-t)t^{-v-1} h_v(t)
$$
where 
$$
h_v(t) \equiv t^{v-1}(t+1) \exp\left( v(1-v) \frac{(t-1)^2}{t} \right) -1 \geq 0.
$$
The last inequality is due to Lemma \ref{lemma_01} in the below.
Therefore we have  $\frac{df_v(t)}{dt} \geq 0$ if $0 < t \leq 1$ and $\frac{df_v(t)}{dt} \leq 0$ if $t \geq 1$. Thus we have $f_v(t) \leq f_v(1) =0$ which implies the inequality (\ref{ineq01_theorem01}).

\hfill \qed

\begin{Lem} \label{lemma_01}
For $t >0$ and $0 \leq v \leq 1$, 
we have
$$
t^{v-1}(t+1) \exp\left( v(1-v) \frac{(t-1)^2}{t} \right)  \geq 1.
$$
\end{Lem}

{\it Proof}:
For any $t >0$ and $0 \leq v \leq 1$, $\exp\left( v(1-v) \frac{(t-1)^2}{t} \right) \geq 1$. In addition, $t^{v-1}(t+1) = t^v + t^{v-1} \geq 1$ since
$t^v \geq 1$, $t^{v-1} >0$ for $t \geq 1$ and $t^v > 0$, $t^{v-1} \geq 1$ for $0<t \leq 1$.
Therefore we have the desired inequality.

\hfill \qed

\begin{remark} \label{remark2_3}
It is known the following inequality (see \cite{ZSF2011,LWZ2015}),
\begin{equation} \label{ineq_ZSF2011_LWZ2015}
 K^r(t) \leq \frac{(1-v)+vt}{t^v} \leq K^R(t)
\end{equation}
where $K(t) = \frac{(t+1)^2}{4t}$ is the Kantorovich constant, $r=\min \{v,1-v\}$ and $R = \max \{v,1-v\}$.
By the numerical computations, we find that there is no ordering between 
 $\exp\left( v(1-v) \frac{(t-1)^2}{t} \right)$ and  $K^R(t)$ so that Theorem \ref{theorem_2_1} is not trivial one.
Actually, we set the function as
$$
l(t,v) \equiv  K^R(t)- \exp\left( v(1-v) \frac{(t-1)^2}{t} \right)
$$
for $t >0$ and $0 \leq v \leq 1$. Then we have $l(1/2,1/5) \simeq 0.0155215$ and $l(1/4,1/5) \simeq -0.00425113$. (In addition, we have similarly $l(3,1/5) \simeq 0.0209862$ and $l(5,1/5) \simeq -0.0682639$.)
\end{remark}

Dragomir also established the following refined Young inequalities with the general inequalities in his paper \cite{Dragomir02}.

\begin{theorem} {\bf (\cite{Dragomir02})} \label{theorem02}
Let $t > 0$ and $0 \leq v \leq 1$.
\begin{itemize}
\item[(i)] If $0 < t \leq 1$, then 
\begin{equation} \label{ineq01_theorem02}
  \exp\left( \frac{v(1-v)}{2}(t-1)^2\right) \leq \frac{(1-v) + v t}{t^v} \leq  \exp\left( \frac{v(1-v)}{2}\left(\frac{1}{t}-1\right)^2\right) 
\end{equation}
\item[(ii)] If $t \geq 1$, then
\begin{equation} \label{ineq02_theorem02}
 \exp\left( \frac{v(1-v)}{2}\left(\frac{1}{t}-1\right)^2\right)  \leq \frac{(1-v) + v t}{t^v} \leq  \exp\left( \frac{v(1-v)}{2}\left(t-1\right)^2\right)
\end{equation}
\end{itemize}
\end{theorem}
%%%%%%%%%%%%%%%%%%%%%%%%%%%%%%%%Proof%%%%%%%%%%%%%%%%%%%%%%%%%%%%
{\it Proof}:
\begin{itemize}
\item[(i)] We prove the first inequality of (\ref{ineq01_theorem02}).
To do this, we set
$$
f_1(t,v) \equiv (1-v)t^{-v}+vt^{1-v}-\exp\left( \frac{v(1-v)(t-1)^2}{2}\right).
$$ 
Then we calculate
$$
\frac{df_1(t,v)}{dt} =v(1-v)(1-t) h_1(t,v),
$$
where
$$
h_1(v,t) \equiv \exp\left( \frac{v(1-v)(t-1)^2}{2}\right)-t^{-v-1}.
$$
From Lemma \ref{lemma_02} in the below, $h_1(t,v) \leq 0$ which means 
$\frac{df_1(t,v)}{dt} \leq 0$. Thus we have $f_1(t,v) \geq f_1(1,v)=0$ which means
the first inequality of (\ref{ineq01_theorem02}) hold for $0< t \leq 1$ and $0 \leq v \leq 1$.

We prove the second inequality of (\ref{ineq01_theorem02}).
To do this, we set
$$
f_2(t,v) \equiv (1-v)t^{-v}+vt^{1-v}-\exp\left( \frac{v(1-v)(t-1)^2}{2t^2}\right).
$$
Then we calculate
$$
\frac{df_2(t,v)}{dt} =v(1-v)(1-t)t^{-3} h_2(t,v),
$$
where
$$
h_2(v,t) \equiv \exp\left( \frac{v(1-v)(t-1)^2}{2t^2}\right)-t^{2-v}.
$$
Since $\exp(u) \geq 1 + u$ for $u \geq 0$, we have $h_2(v,t) \geq g_2(v,t)$, where
$$
g_2(t,v) \equiv 1+\frac{v(1-v)(t-1)^2}{2t^2} -t^{2-v}.
$$
Then we calculate
$$
\frac{dg_2(t,v)}{dt} =t^{-3}\left\{ (v-2) t^{4-v} +v(1-v)(t-1) \right\} \leq 0
$$
for $0< t \leq 1$ and $0 \leq v \leq 1$.
Thus we have $g_2(t,v) \geq g_2(1,v) =0$ which implies $h_2(t,v) \geq 0$ which means
$\frac{df_2(t,v)}{dt} \geq 0$. Thus we have $f_2(t,v) \leq f_2(1,v) =0$ which means
the second inequality of (\ref{ineq01_theorem02}) hold for $0< t \leq 1$ and $0 \leq v \leq 1$.
\item[(ii)] We prove the first inequality of (\ref{ineq02_theorem02}).
The functions $f_2(t,v)$, $h_2(t,v)$ and $g_2(t,v)$ were defined in the process of the proof of the second inequality in (i).
Here we set the function for $t \geq 1$
$$
k(t,v) \equiv (v-2)t^{4-v}+v(1-v)(t-1).
$$
Then we calculate
$$
\frac{dk(t,v)}{dt} =(v-2)(4-v)t^{3-v}+v(1-v),\,\,\frac{d^2k(t,v)}{dt^2} =(v-2)(4-v)(3-v)t^{2-v} \leq 0.
$$
Thus we have $\frac{dk(t,v)}{dt} \leq \frac{dk(1,v)}{dt} =-2\left(v-\frac{7}{4} \right)^2 -\frac{15}{8} \leq 0$ so that we have $k(t,v) \leq k(1,v) = v-2 <0$.
Thus we have $\frac{dg_2(t,v)}{dt} \leq 0$ for $t \geq 1$ so that $g_2(t,v) g_2(1,v)=0$ which implies $h_2(t,v) \geq 0$. Therefore we have $\frac{df_2(t,v)}{dt} \leq 0$ so that
$f_2(t,v) \geq f_2(1,v) =0$ which implies the first inequality of (\ref{ineq02_theorem02}).

We prove the second inequality of (\ref{ineq02_theorem02}).
The functions $f_1(t,v)$, $h_1(t,v)$ and $g_1(t,v)$ were defined in the process of the proof of the first inequality in (i).
For $t \geq 1$, we easily find $\frac{dg_1(t,v)}{dt} =v(1-v)(t-1) +(v+1)t^{v-2} \geq 0$ so that $g_1(t,v) \geq g_1(1,v) =0$. Thus we have $h_1(t,v) \geq 0$ which means $\frac{df_1(t,v)}{dt} \leq 0$. Therefore we have $f(t,v) \leq f(1,v) =0$ 
which implies the second inequality of (\ref{ineq02_theorem02}).
\end{itemize}

\hfill \qed

\begin{Lem} \label{lemma_02}
For $0 < t \leq 1$ and $0 \leq v \leq 1$, we  have
\begin{equation} \label{ineq01_lemma02}
t^{v+1}\exp \left( \frac{v(1-v)}{2} (t-1)^2\right) \leq 1.
\end{equation}
\end{Lem}

{\it Proof}:
We set the function as
$$
f_v(t) \equiv (v+1) \log t +\frac{v(1-v)}{2} (t-1)^2
$$
Then we calculate
$$
\frac{df_v(t)}{dt} =\frac{v+1}{t}+v(1-v)(t-1),\,\,\frac{d^2f(t,v)}{dt^2}=-\frac{v+1}{t^2} +v(1-v),\,\,\frac{d^3f(t,v)}{dt^3} =\frac{2(v+1)}{t^3} \geq 0
$$
Thus we have $\frac{d^2f_v(t)}{dt^2} \leq \frac{d^2f_v(1)}{dt^2} =-v^2-1 \leq 0$ so that we have $\frac{df_v(t)}{dt} \geq \frac{df_v(1)}{dt} =v+1 \geq 0$. Therefore we have $f_v(t) \leq f_v(1) =0$ which implies the inequality (\ref{ineq01_lemma02}).

\hfill \qed

\begin{remark}
The second inequalities (\ref{ineq01_theorem02}) and  (\ref{ineq02_theorem02})  refine the second inequality in \cite[Corollary2.2 (i)]{FM2011}.
\end{remark}

We obtained the following results in our previous paper.
\begin{proposition} {\bf (\cite{FM})} \label{proposition_FM}
For $0 \leq v \leq 1$ and $0< t \leq 1$, we have
$$
m_v(t)  \leq  \frac{(1-v) + vt}{t^v} \leq  M_v(t),
$$
where
$$
m_v(t) \equiv 1+\frac{v(1-v)(t-1)^2}{2}\left(\frac{t+1}{2}\right)^{-v-1}, \quad M_v(t) \equiv 1+\frac{v(1-v)(t-1)^2}{2}t^{-v-1}.
$$
\end{proposition}

\begin{remark}
As shown in our previous paper \cite{FM}, we have the inequality
$$
M_v(t) \leq \exp \left(v(1-v)\frac{(t-1)^2}{t} \right)
$$
for $0<t \leq 1$ and $0 \leq v \leq \frac{1}{2}$. That is, the second inequality in Proposition \ref{proposition_FM} gives better bound than the inequality (\ref{ineq01_theorem01}), in case of $0<t \leq 1$ and $0 \leq v \leq \frac{1}{2}$.
\end{remark}

In the following proposition, we give the comparison on bounds in (i) Theorem \ref{theorem02} and in Proposition \ref{proposition_FM}.

\begin{proposition} \label{proposition2_8}
For $0 \leq v \leq 1$ and $0< t \leq 1$, we have
\begin{equation} \label{ineq01_theorem_2_8}
 M_v(t) \leq \exp\left( \frac{v(1-v)}{2} \left( \frac{1}{t} -1\right)^2 \right)
\end{equation}
and
\begin{equation} \label{ineq02_theorem_2_8}
\exp\left( \frac{v(1-v)}{2} \left( t -1\right)^2 \right) \leq m_v(t)
\end{equation}
\end{proposition}

{\it Proof}:
We use the inequality 
$$
\exp x \geq 1+x +\frac{1}{2} x^2,\quad (x \geq 0).
$$
Then we calculate
\begin{eqnarray*}
&&\exp\left( \frac{v(1-v)}{2} \left( \frac{1}{t} -1\right)^2 \right) -1-\frac{v(1-v)(t-1)^2}{2t^{v+1}}\\
&&\geq \frac{v(1-v)(t-1)^2}{2t^2} \left( \frac{t^{v-1}-1}{t^{v-1}} +\frac{v(1-v)(t-1)^2}{4t^2}\right) \geq 0
\end{eqnarray*}
for $0 \leq v \leq 1$ and $0< t \leq 1$. Thus the inequality (\ref{ineq01_theorem_2_8}) was proved.

Putting $s = \frac{t+1}{2}$, the inequality (\ref{ineq02_theorem_2_8})  is equivalent to the inequality
\begin{equation} \label{ineq01_proof_theorem_2_8}
\exp\left( 2v(1-v)(s-1)^2\right) \leq 1+2v(1-v)(s-1)^2s^{-v-1},\quad \left(\frac{1}{2} < s \leq 1,\,\, 0 \leq v \leq 1 \right).
\end{equation}
For the special case $v=0,1$ or $s=1$, the equality holds in (\ref{ineq01_proof_theorem_2_8}) so that we assume $2v(1-v)(s-1)^2 \ne 0$.
Then we use the inequality
$$
\exp x < \frac{1}{1-x},\quad (x <1).
$$
We calculate
\begin{eqnarray*}
&&m_v(t) - \exp\left( \frac{v(1-v)}{2} \left( t -1\right)^2 \right)  = 1+2v(1-v)(s-1)^2s^{-v-1} - \exp\left( 2v(1-v)(s-1)^2\right)\\
&& > 1+2v(1-v)(s-1)^2s^{-v-1} -\frac{1}{1-2v(1-v)(s-1)^2} 
 = \frac{2v(1-v)(s-1)^2s^{-v-1} g_v(s)}{1-2v(1-v)(s-1)^2},
\end{eqnarray*}
where 
$$
g_v(s) \equiv 1-s^{v+1} -2v(1-v)(s-1)^2.
$$
We prove $g_v(s) \geq 0$. To this end, we calculate
\begin{eqnarray*}
&& g_v'(s)=-(v+1)s^v-4v(1-v)(s-1),\quad  g''_v(s)=-v(v+1)s^{v-1} -4v(1-v) \\
&& g^{(3)}_v(s) =v(1-v)(v+1)s^{v-2} \geq 0.
\end{eqnarray*}
Thus we have $g''_v(s) \leq g''_v(1) = v(3v-5) \leq 0 $ so that we have
$g'_v(s) \leq g'_v(1/2) =2 h(v)$, where $h(v) \equiv v(1-v)-\frac{v+1}{2^{v+1}}$.
From Lemma \ref{lemma03}, $h(v) \leq 0$ for $0 \leq v \leq 1$.
Thus we have $g'_v(s) \leq 0$ which implies $g_v(s) \geq g_v(1) =0$.
Therefore we have
$$
m_v(t) - \exp\left( \frac{v(1-v)}{2} \left( t -1\right)^2 \right) >0
$$ 
for $0<v<1$ and $0<t<1$. Taking account for the equality cases happen if $v=0,1$ or $t=1$, we have the inequality (\ref{ineq02_theorem_2_8}).

\hfill \qed 

\begin{Lem}\label{lemma03}
For $0 \leq v \leq 1$, we have 
$$
\frac{v+1}{2^{v+1}} \geq v(1-v).
$$
\end{Lem}

{\it Proof}:
Since $v(1-v) \leq \frac{1}{4}$ for  $0 \leq v \leq 1$, it is sufficient to prove
$\frac{v+1}{2^{v+1}}  \geq \frac{1}{4}$. So we put $l(v) \equiv 2(v+1)-2^v$.
Then we have $l''(v) =-(\log 2)^2 2^{v} \leq 0$, $l(0)=1$ and $l(1)=2$.
Therefore we have $l(v) \geq 0$.

\hfill \qed 

\begin{remark}
As for the bounds on the ratio of arithmetic mean to geometric mean $\frac{(1-v)+v t}{t^v}$,
Proposition \ref{proposition2_8} shows Proposition \ref{proposition_FM} is better than (i) of Theorem \ref{theorem02}, for the case $0<t \leq 1$.
\end{remark}
%%%%%%%%%%%%%%%%%%%%%%%%%%%%%%%%%%%%%%%%%%%%%%%%%%%%%%%%%%%%%
%%%%%%%%%%%%%%%%%%%%%%%%%%%%%%%%%%%%%%%%%%%%%%%%%%%%%%%%%%%%%
%%%%%%%%%%%%%%%%%%%%%%%%%%%%%%%%%%%%%%%%%%%%%%%%%%%%%%%%%%%%%
%%%%%%%%%%%%%%%%%%%%%%%%%%%%%%%%%%%%%%%%%%%%%%%%%%%%%%%%%%%%%
%%%%%%%%%%%%%%%%%%%%%%%%%%%%%%%%%%%%%%%%%%%%%%%%%%%%%%%%%%%%%
%%%%%%%%%%%%%%%%%%%%%%%%%%%%%%%%%%%%%%%%%%%%%%%%%%%%%%%%%%%%%

\section{Further improvement of Young inequality}

We give new improvement of Young inequality which is a further improvement of Theorem \ref{theorem_2_1}. Throughout this section, we use the generalized exponential function defined by $\exp_r(x) \equiv (1+r x)^{1/r}$ for $x >0$ and $-1 \leq r \leq 1$ with $r \ne 0$ under the assumption that $1+r x \geq 0$.

\begin{theorem} \label{theorem3_1}
For $t >0$ and $0 \leq v \leq 1$, 
\begin{equation} \label{ineq01_theorem3_1}
\frac{(1-v) + vt}{t^v} \leq 1+ v(1-v) \frac{(t-1)^2}{t} 
\end{equation}
\end{theorem}

{\it Proof}:
We set the function 
$$
f_v(t) \equiv 1+ v(1-v) \frac{(t-1)^2}{t} -\frac{(1-v) + vt}{t^v}
$$
for $t >0$ and $0 \leq v \leq 1$. Then we have
$$
\frac{df_v(t)}{dt} = \frac{v(1-v)(t-1)}{t^{v+1}} \left(t^v +t^{v-1} -1 \right).
$$
Since $t^v \geq 1$ for $t \geq 1$ and $t^{v-1} \geq 1$ for $0 < t \leq 1$, $t^v + t^{v-1} -1 \geq 0$. Thus we have $\frac{df_v(t)}{dt} = 0$ when $t =1$  and
$\frac{df_v(t)}{dt} \leq 0$ for $0<t<1$, and $\frac{df_v(t)}{dt} \geq 0$ for $t>1$.
Therefore we have $f_v(t) \geq f_v(1) =0$.

\hfill \qed 

\begin{Lem} \label{lemma3_2}
The function $\exp_r(x)$ defined for $x >0$ and $0 < r \leq 1$, is monotone decreasing in $r.$
\end{Lem}
{\it Proof}:
We calculate 
$
\frac{d \exp_r(x)}{dx} = \frac{(1+r x)^{\frac{1-r}{r}} g(rx)}{r^2}
$
where
$g(y) \equiv y -(1+y)\log (1 + y)$ for $y >0$. 
Since $\frac{dg(y)}{dy} = -\log (1+y) <0$, $f(y) \leq f(0) =0$. We thus have
$\frac{d \exp_r(x)}{dx} \leq 0$.

\hfill \qed 

\begin{corollary} \label{cor3_3}
For $t >0$, $0 \leq v \leq 1$ and $0<r \leq 1$, 
\begin{equation} \label{ineq01_corollary3_3}
\frac{(1-v) + vt}{t^v} \leq \exp_r\left( v(1-v)\frac{(t-1)^2}{t}\right).
\end{equation}
\end{corollary}

{\it Proof}:
Since $\exp_1\left( v(1-v)\frac{(t-1)^2}{t}\right)= 1+ v(1-v)\frac{(t-1)^2}{t}$, we have the desired result by Theorem \ref{theorem3_1} and Lemma \ref{lemma3_2}.

\hfill \qed

\begin{remark} 
Since $\lim_{r\to 0} \exp_r(x) = \exp(x)$ and Lemma \ref{lemma3_2}, we have
$$
1+ v(1-v)\frac{(t-1)^2}{t} \leq \exp \left( v(1-v)\frac{(t-1)^2}{t}\right)
$$
which means that the right hand side in Theorem \ref{theorem3_1} gives the tighter upper bounds of $\frac{(1-v)+vt}{t^v}$ than one in Theorem \ref{theorem_2_1}.
\end{remark}

\begin{remark}
Proposition \ref{proposition_FM} shows the upper bound of $\frac{(1-v)+vt}{t^v}$ is $M_v(t)$ for $0 < t \leq 1$, while Theorem \ref{theorem3_1} gives the upper bound of $\frac{(1-v)+vt}{t^v}$ for all $t >0$. In addition, for the case $t^v \leq \frac{1}{2}$, the right hand side in Theorem \ref{theorem3_1} gives the tighter upper bounds of $\frac{(1-v)+vt}{t^v}$ than $M_v(t)$ in Proposition \ref{proposition_FM}. 
\end{remark}

\begin{remark}
It is easy to see that $1+v(1-v)\frac{(t-1)^2}{t} \leq 1+\frac{(t-1)^2}{4t} = K(t)$.
As we noted in Remark \ref{remark2_3}, the inequalities (\ref{ineq_ZSF2011_LWZ2015}) are known. By the numerical computations, we have no ordering between $K^R(t)$ and $1+v(1-v)\frac{(t-1)^2}{t}$. Actually, we set the function for $t>0$ and $0\leq v \leq 1$ as
$$
u_1(t,v) \equiv K^R(t)-\left(1+v(1-v)\frac{(t-1)^2}{t}\right).
$$
Then $u_1(1/2,0.6) \simeq -0.0467732$ and $u_1(1/2,0.9) \simeq 0.0668271$.
\end{remark}

\begin{remark}
From the proof of Lemma \ref{lemma3_2}, we find that the function $\exp_r(x)$ is monotone decreasing for $r>0$. However the following inequality does not hold in general
$$
\frac{(1-v)+vt}{t^v} \leq \exp_r\left( v(1-v)\frac{(t-1)^2}{t}\right),\quad (r>1)
$$
since we have a counter-example. For example, if we take $v=0.999999$ and $t=0.000001$, then $\exp_{1.001}\left( v(1-v)\frac{(t-1)^2}{t}\right)-\frac{(1-v)+vt}{t^v} \simeq -0.000360488$. This example suggests the optimality of $r$ satisfying the inequality (\ref{ineq01_corollary3_3}) is equal to $1$. 
\end{remark}

%This remark suggests the optimality of $r$ satisfying the inequality %(\ref{ineq01_corollary3_3}) is equal to $1$. That is, we have the following %proposition.
%\begin{proposition}
%\begin{equation}
%\sup\left\{ r >0 \left|  \frac{(1-v) + vt}{t^v} \leq \exp_r\left( v(1-v)\frac{(t-1)^2}{t}\right)\right. \right\} =1.
%\end{equation}
%\end{proposition} 

%{\it Proof}:
%For arbitrary $\varepsilon >0$, we consider the inequality
%\begin{eqnarray*}
%&&\frac{(1-v) + vt}{t^v} \leq \exp_{1+\varepsilon}\left( v(1-v)\frac{(t-1)^2}{t}\right)
%=\left( 1+(1+\varepsilon)v(1-v)\frac{(t-1)^2}{t}\right)^{\frac{1}{1+\varepsilon}}\\
%&& \leq \left( 1+(1+\varepsilon)v(1-v)\frac{(t-1)^2}{t}\right),
%\end{eqnarray*}
%which implies
%$$\varepsilon \geq \frac{t}{v(1-v)(t-1)^2} \left\{\frac{(1-v)+vt}{t^v} -1-v(1-v)\frac{(t-1)^2}{t}\right\}.$$
%When $v \to  1$, the right hand side in the above inequality converges to %$\frac{t(1-t+\log t)}{(t-1)^2}$ which converges to $1$ in the limit $t \to 0$.
%The we obtain $\varepsilon \geq 0$ which implies the proposition.

%\hfill \qed

As similar way to the above, we can improve Theorem \ref{theorem02} in the following.

\begin{theorem} \label{theorem3_6}
Let $t > 0$ and $0 \leq v \leq 1$.
\begin{itemize}
\item[(i)] If $0 < t \leq 1$, then
$$
\frac{1}{1-\frac{v(1-v)}{2}(t-1)^2} \leq \frac{(1-v)+vt}{t^v} \leq 1+ \frac{v(1-v)}{2}\left(\frac{1}{t} -1\right)^2
$$
\item[(ii)] If $t \geq 1$, then
$$
\frac{1}{1-\frac{v(1-v)}{2}(\frac{1}{t}-1)^2} \leq \frac{(1-v)+vt}{t^v} \leq 1+ \frac{v(1-v)}{2}\left(t -1\right)^2
$$
\end{itemize}

{\it Proof}:
Firstly, we prove the second inequality in (i). To this end, we set the function as
$$
f_v(t) \equiv 1+\frac{v(1-v)}{2}\left(\frac{1}{t} -1\right)^2-\frac{(1-v)+vt}{t^v}
$$
for $0 < t \leq 1$ and   $0 \leq v \leq 1$.
Then we find by elementary calculation
$$
\frac{df_v(t)}{dt} = v(1-v)(t-1)t^{v+3} (t^v -t^2) \leq 0
$$
so that $f_v(t) \geq f_v(1)=0$.
Secondary we prove the first inequality in (i). To this end, we set the function as
$$
g_v(t) \equiv \left\{(1-v)+vt\right\}\left\{1-\frac{v(1-v)}{2}(t-1)^2\right\}-t^v
$$
for $0 < t \leq 1$ and   $0 \leq v \leq 1$.
Then we find by elementary calculations
\begin{eqnarray*}
&& \frac{dg_v(t)}{dt} = v \left\{-2t^{v-1}-3v(1-v)t^2 -2(3v^2-4v+1)t+(3v^2-5v+4) \right\}  \\
&& \frac{d^2g_v(t)}{dt^2} = v(1-v) \left\{ t^{v-2}-3vt +(3v-1) \right\},\quad
 \frac{d^3g_v(t)}{dt^3} =v(1-v) \left\{ (v-2)t^{v-3} -3v \right\} \leq 0.
\end{eqnarray*}
Thus we have $\frac{d^2g_v(t)}{dt^2} \geq \frac{d^2g_v(1)}{dt^2}=0$ so that
$\frac{dg_v(t)}{dt} \leq \frac{dg_v(1)}{dt}=0$ which implies $f_v(t) \geq f_v(1) =0$.

Thirdly we prove the first inequality in (ii). To this end, we set the function as
$$
h_v(t) \equiv \left\{(1-v)+vt\right\}\left\{1-\frac{v(1-v)}{2}\left(\frac{1}{t}-1\right)^2\right\}-t^v
$$
for $t \geq 1$ and   $0 \leq v \leq 1$.
Then we find by elementary calculations
\begin{eqnarray*}
&& \frac{dh_v(t)}{dt} =\frac{v}{2t^3} \left\{-2t^{v+2}+(v^2-v+2)t^3-(3v^2-5v+2)t+2(1-v)^2 \right\}, \\
&& \frac{d^2h_v(t)}{dt^2} = \frac{v(1-v)}{t^4} l_v(t),\quad l_v(t)\equiv t^{v+2}+(2-3v)t+3v-3.
\end{eqnarray*}
Since $\frac{dl_v(t)}{dt} =(v+2)t^{v+1} +2-3v$ and $\frac{d^2l_v(t)}{dt^2} = (v+2)(v+1)t^v \geq 0$, we have $\frac{dl_v(t)}{dt} \geq \frac{dl_v(1)}{dt} =4-2v \geq 0$ so that $l_v(t) \geq l_v(1)=0$ which implies $\frac{d^2h_v(t)}{dt^2}  \geq 0$.
Thus we have $\frac{dh_v(t)}{dt} \geq \frac{dh_v(1)}{dt} =0$ so that $h_v(t) \geq h_v(1) =0$.

Finally we prove the second inequality in (ii). To this end, we set the function as
$$
k_v(t) \equiv 1+\frac{v(1-v)}{2}\left(t -1\right)^2-\frac{(1-v)+vt}{t^v}
$$
for $t \geq 1$ and   $0 \leq v \leq 1$.
Since $\frac{dk_v(t)}{dt} =v(1-v)(t-1)(1-t^{-v-1}) \geq 0$, we have $k_v(t) \geq k_v(t) =0$.

\hfill \qed

\end{theorem}

\begin{Lem} \label{lemma3_7}
The function $\exp_r(x)$ defined for $0\leq x \leq 1$ and $-1 \leq r <0$, is monotone decreasing in $r.$
\end{Lem}

{\it Proof}:
We calculate 
$
\frac{d \exp_r(x)}{dx} = \frac{(1+r x)^{\frac{1-r}{r}} g(rx)}{r^2}
$
where $g(y) \equiv y -(1+y)\log (1 + y)$ for $-1 \leq y \leq 0$. 
Since $\frac{dg(y)}{dy} = -\log (1+y) >0$, $f(y) \leq f(0) =0$. We thus have
$\frac{d \exp_r(x)}{dx} \leq 0$.

\hfill \qed

\begin{corollary} \label{cor3_8}
Let $t >0$, $0 \leq v \leq 1$, $-1 \leq r_1 <0$ and $0< r_2 \leq 1$.
\begin{itemize}
\item[(i)] If $0 < t \leq 1$, then
$$
\exp_{r_1}\left(\frac{v(1-v)}{2}(t-1)^2 \right) \leq \frac{(1-v)+vt}{t^v} \leq \exp_{r_2}\left(\frac{v(1-v)}{2}\left(\frac{1}{t}-1\right)^2 \right).
$$
\item[(ii)] If $t \geq 1$, then
$$
\exp_{r_1}\left(\frac{v(1-v)}{2}\left(\frac{1}{t}-1\right)^2 \right)\leq \frac{(1-v)+vt}{t^v} \leq
\exp_{r_2}\left(\frac{v(1-v)}{2}(t-1)^2 \right).
$$
\end{itemize} 
\end{corollary}

{\it Proof}:
Taking account for $\exp_1(x) = 1+x$, applying two second inequalities in (i) and (ii) of Theorem \ref{theorem3_6} and Lemma \ref{lemma3_2}, we obtain two second inequalities in (i) and (ii) of this theorem.

Taking account for $\exp_{-1}(x) =\frac{1}{1-x}$ and $\frac{v(1-v)}{2}(t-1)^2 <1$ for $0< t \leq 1$, $\frac{v(1-v)}{2}\left(\frac{1}{t} -1 \right)^2 <1$ for $t \geq 1$,
applying two first inequalities in (i) and (ii) of Theorem \ref{theorem3_6} and Lemma \ref{lemma3_7}, we obtain two first inequalities in (i) and (ii) of this theorem.

\hfill \qed

\begin{remark}
As we noted $\lim_{r \to 0} \exp_r (x) =\exp(x)$, and Lemma \ref{lemma3_2} and Lemma \ref{lemma3_7} assure that Theorem \ref{theorem3_6} gives tighter bounds of $\frac{(1-v)+vt}{t^v}$ than Theorem \ref{theorem02}.
\end{remark}

\begin{remark} \label{remark3_10}
Bounds in (i) of Theorem \ref{theorem3_6} can be compared with Proposition \ref{proposition_FM}. As for upper bound, $M_v(t)$ gives tighter than the right hand side in the second inequality of (i) of  Theorem \ref{theorem3_6}.
Since $\left(\frac{t+1}{2}\right)^{v+1} \geq t^2$ for $0<t \leq 1$, $m_v(t)$ also gives tighter than the left hand side in the first inequality of (i) of  Theorem \ref{theorem3_6}. The inequality $\left(\frac{t+1}{2}\right)^{v+1} \geq t^2$ can be proven in the following. We set $f_v(t)=(v+1)\log\frac{t+1}{2} -2\log t$, then we have $\frac{df_v(t)}{dt} = \frac{(v-1)t-2}{t(t+1)} \leq 0$ so that $f_v(t) \geq f_v(1) =0$. 
\end{remark}

Note that Remark \ref{remark3_10} and (i) of Corollary \ref{cor3_8} give Proposition \ref{proposition2_8}. However, Corollary \ref{cor3_8} gives alternative tight bounds of $\frac{(1-v)+vt}{t^v}$ when $t \geq 1$.

\begin{remark}
We give comparisons our inequalities obtained in Theorem \ref{theorem3_6} with the inequalities (\ref{ineq_ZSF2011_LWZ2015}). By the numerical computations, we have no ordering between $K^R(t)$ and $1+v(1-v)\frac{(t-1)^2}{2t^2}$ for $0<t\leq 1$. Actually, we set the function for $0<t\leq 1$ and $0\leq v \leq 1$ as
$$
u_2(t,v) \equiv K^R(t)-\left(1+v(1-v)\frac{(t-1)^2}{2t^2}\right).
$$
Then $u_2(1/2,0.6) \simeq -0.0467732$ and $u_2(1/2,0.9) \simeq 0.0668271$.
(We easily find that $u_1(1/2,v)=u_2(1/2,v)$.)
We also have no ordering between $K^R(t)$ and $1+v(1-v)\frac{(t-1)^2}{2}$ for $t \geq 1$. Actually, we set the function for $t \geq 1$ and $0\leq v \leq 1$ as
$$
u_3(t,v) \equiv K^R(t)-\left(1+v(1-v)\frac{(t-1)^2}{2}\right).
$$
Then $u_3(2,0.6) \simeq -0.0467732$ and $u_3(2,0.9) \simeq 0.0668271$.
(We easily find that $u_2(1/2,v)=u_1(1/2,v)=u_1(2,v)=u_3(2,v)$.)
As for the lower bounds, we have no ordering between $K^r(t)$ and $\left(1-\frac{v(1-v)}{2}(t-1)^2\right)^{-1}$ for $0<t\leq 1$. Actually, we set the function for $0<t\leq 1$ and $0\leq v \leq 1$ as
$$
l_1(t,v) \equiv K^r(t)-\left(1-\frac{v(1-v)}{2}(t-1)^2\right)^{-1}.
$$
Then $l_1(3/5,0.1) \simeq -0.000777493$ and $l_1(3/5,0.4) \simeq 0.00657566$.
We also have no ordering between $K^r(t)$ and $\left(1-\frac{v(1-v)}{2t^2}(t-1)^2\right)^{-1}$ for $t \geq 1$. Actually, we set the function for $t  \geq 1$ and $0\leq v \leq 1$ as
$$
l_2(t,v) \equiv K^r(t)-\left(1-\frac{v(1-v)}{2t^2}(t-1)^2\right)^{-1}.
$$
Then $l_2(5/3,0.1) \simeq -0.000777493$ and $l_2(5/3,0.4) \simeq 0.00657566$.
\end{remark}

We close this section showing operator versions for Corollary \ref{cor3_3} and Corollary \ref{cor3_8} in the following. We denote weighted arithmetic mean and geometric mean for two strictly positive operators $A$ and $B$ by
$A \nabla_v B$ and $A\#_v B$, respectively.

\begin{corollary}
Let $0\leq v \leq 1$, $0< r \leq 1$ and let $A$ and $B$ be strictly positive operators satisfying (i) $0<m \leq A \leq m' < M' \leq B \leq M$ or (ii) $0 < m \leq B \leq m' < M' \leq A \leq M$ with $h =\frac{M}{m}$ and $h'=\frac{M'}{m'}$.
Then
$$
A \nabla_v B \leq \exp_r\left(4v(1-v)(K(h)-1) \right) A\#_vB,
$$
where $K(h) \equiv \frac{(h+1)^2}{4h}$ is Kantrovich constant.
\end{corollary}
 
{\it Proof}:
The inequality (\ref{ineq01_corollary3_3}) is equivalent to
$$
(1-v) + vt \leq \exp_r\left(4v(1-v)(K(t)-1)\right) t^v
$$ 
for any $t >0$. Thus we have the following operator inequality for strictly positive operator $T$ such that $0<m \leq T \leq M$,
$$
(1-v) + vT \leq \max_{m\leq t \leq M}\exp_r\left(4v(1-v)(K(t)-1)\right) T^v.
$$
Here we put $T = A^{-1/2}BA^{-1/2}$.
In the case of (i),  we have $h' \leq A^{-1/2}BA^{-1/2} \leq h$. Then we have
$$
(1-v) + vA^{-1/2}BA^{-1/2} \leq \max_{h'\leq t \leq h}\exp_r\left(4v(1-v)(K(t)-1)\right) (A^{-1/2}BA^{-1/2})^v.
$$
In the case of (ii), we also have $\frac{1}{h} \leq A^{-1/2}BA^{-1/2} \leq \frac{1}{h'}$. Then we also have
$$
(1-v) + vA^{-1/2}BA^{-1/2} \leq \max_{1/h\leq t \leq 1/h'}\exp_r\left(4v(1-v)(K(t)-1)\right) (A^{-1/2}BA^{-1/2})^v.
$$
Since $K(h)$ is decreasing for $0< h \leq 1$, increasing for $h \geq 1$ and $K(1/h)=K(h) \geq K(1)=1$, we obtain the desired result by multiplying $A^{1/2}$ to both sides in two above inequalities.  

\hfill \qed

\begin{corollary}
Let $0\leq v \leq 1$, $-1\leq r_1 <0$, $0< r_2 \leq 1$ and let $A$ and $B$ be strictly positive operators satisfying (i) $0<m \leq A \leq m' < M' \leq B \leq M$ or (ii) $0 < m \leq B \leq m' < M' \leq A \leq M$ with $h =\frac{M}{m}$ and $h'=\frac{M'}{m'}$.
Then
$$
\exp_{r_1}\left( \frac{v(1-v)}{2} \left(\frac{h-1}{h}\right)^2 \right)A\#_vB \leq A\nabla_v B \leq \exp_{r_2}\left( \frac{v(1-v)}{2} \left(h'-1\right)^2 \right)A\#_vB 
$$
\end{corollary}

{\it Proof}:
The inequalities in (i) and (ii) of Corollary \ref{cor3_8} can be written as
$$
\exp_{r_1}\left( \frac{v(1-v)}{2} \left(1-\frac{\min\{1,t\}}{\max\{1,t\}}\right)^2 \right) \leq \frac{(1-v)+vt}{t^v} \leq \exp_{r_2}\left( \frac{v(1-v)}{2} \left(1-\frac{\max\{1,t\}}{\min\{1,t\}}\right)^2 \right). 
$$
The rest of the proof goes similar way to the proof of \cite[Corollary 1]{Dragomir02}. We omit its details.

\hfill \qed

\section*{Acknowledgement}
The author was partially supported by JSPS KAKENHI Grant Number
16K05257.

\end{document}